\documentclass[12pt,a4paper,english,intoc,bibliography=totoc,index=totoc,BCOR10mm,captions=tableheading,titlepage,fleqn]{scrbook}
\usepackage{lmodern}

\usepackage[T1]{fontenc}
\usepackage[latin9]{inputenc}
\usepackage{babel}
\usepackage{amsmath}
\usepackage{amssymb}
\usepackage[unicode=true,
 bookmarks=true,bookmarksnumbered=true,bookmarksopen=true,bookmarksopenlevel=2,
 breaklinks=false,pdfborder={0 0 0},backref=false,colorlinks=false]
 {hyperref}
\hypersetup{pdftitle={Abstract},
 pdfauthor={Uwe Stöhr},
 pdfpagelayout=OneColumn, pdfnewwindow=true, pdfstartview=XYZ, plainpages=false}
\usepackage{breakurl}

\makeatletter

\usepackage{calc}

\makeatother

\begin{document}
\noindent \begin{center}
\textbf{\Large On the definability of mad families of vector spaces}
\par\end{center}{\Large \par}

\noindent \begin{center}
{\large Haim Horowitz and Saharon Shelah}
\par\end{center}{\large \par}

\noindent \begin{center}
\textbf{\small Abstract}
\par\end{center}{\small \par}

\noindent \begin{center}
{\small We consider the definability of mad families in vector spaces
of the form $\underset{n<\omega}{\bigoplus} F$ where $F$ is a field
of cardinality $\leq \aleph_0$. We show that there is no analytic
mad family of subspaces when $F=\mathbb{F}_2$, partially answering
a question of Smythe. Our proof relies on a variant of Mathias forcing
restricted to a certain idempotent ultrafilter whose existence follows
from Glazer's proof of Hindman's theorem.}%
\footnote{{\small Date: November 6, 2018}{\small \par}



Publication 1151 of the second author%
}
\par\end{center}{\small \par}

Assuming the axiom of choice, one can easily construct sets of reals
exhibiting certain maximality properties. One classical example of
such sets is provided by maximal almost disjoint (mad) families. Recall
that $\mathcal A \subseteq [\omega]^{\omega}$ is mad if $x\neq y \in \mathcal A \rightarrow |x\cap y|<\aleph_0$
and $\mathcal A$ is maximal with respect to this property. The study
of the definability of mad families goes back to Mathias in the 1970s.
As with other regularity properties, such as Lebesgue measurability
and the Baire property, it turned out that mad families can't be too
nicely definable:

\textbf{Theorem ({[}Ma{]}): }There are no analytic mad families.

The study of the definability of relatives of mad families has attracted
significant attention in recent years. In a surprising development,
contrary to the pattern described above, it was shown in {[}HwSh:1089{]}
that there exists a Borel maximal eventually different family (where
$\mathcal F \subseteq \omega^{\omega}$ is maximal eventually different
if $f\neq g \in \mathcal F \rightarrow f(n) \neq g(n)$ for large
enough $n$, and $\mathcal F$ is maximal with respect to this property).
In a subsequent work ({[}HwSh:1095{]}), the definability of another
type of relatives of mad families - known as maximal cofinitary groups
- was studied, and it was shown that there exists a Borel maximal
cofinitary group.

The current paper studies the definability of a new variant of mad
families recently introduced by Iian Smythe in {[}Sm{]}. Given an
$\aleph_0$-dimensional vector space $V=\underset{n<\omega}{\bigoplus} F$
over a field $F$ of cardinality $\leq \aleph_0$, we can regard $2^V$
as a Polish space and consider the definability of families of subsets
of $V$. Mad families of subspaces of $V$ will be defined in the
natural way, see Definition 1 below. Our main goal is to provide a
partial answer to the following question:

\textbf{Question ({[}Sm{]}): }For $V$ as above, is there an analytic
mad family of subspaces of $V$?

We shall prove that for $F=\mathbb{F}_2$, the answer is negative,
i.e. we have an analog of Mathias' theorem. We shall assume towards
contradiction that $\mathcal A$ is an analytic mad family of subspaces
of $V$. A main ingredient in our proof will be the existence of a
nonprincipal ultrafilter $D$ on $V$ that is disjoint to $\mathcal A$,
contains all subspaces of finite codimension and has the property
that if $A\in D$, then $v+A \in D$ for $D$-almost all $v$. Such
an ultrafilter will be provided using Glazer's argument for the existence
of idempotent ultrafilters in $\beta(V)$. We shall then consider
the forcing $\mathbb{Q}_D$, a variant of Mathias forcing restricted
to the ultrafilter $D$. $\mathbb{Q}_D$ will introduce a generic
subset $\{ \underset{\sim}{y_k} : k<\omega\}$ of $V$ whose span
is almost contained in every element of $D$. The above invriance
property of $D$ will be used to show that $\{ \underset{\sim}{y_k} : k<\omega\}$
is infinite using a standard density argument. As $\mathcal A$ is
analytic, it remains mad in $\bold{V}^{\mathbb{Q}_D}$, and we can
find a name $\underset{\sim}{B}$ of a new element of $\mathcal A$
that has an infinite intersection with $span \{ \underset{\sim}{y_k} : k<\omega \}$.
We shall then work over a countable elementary submodel $N$ of $H((2^{\aleph_0})^+)$
and construct two generic sets $G_1,G_2$ over $N$ that decide $\underset{\sim}{B}$
in two different ways but still give an infinite intersection of the
two versions of $\underset{\sim}{B}$. By the absoluteness of $\mathcal A$,
this will contradict its almost disjointness. 

\textbf{Definition 1: }a. Let $V$ be an $\aleph_0$-dimensional vector
space over a field $F$ of cardinality $\leq \aleph_0$. We say that
the subspaces $S_1,S_2 \subseteq V$ are almost disjoint if $dim(S_1 \cap S_2)<\aleph_0$.

b. We say that $\mathcal A$ is a mad family of subspaces of $V$
(or a $V$-mad family) if $\mathcal A$ is infinite, the members of
$\mathcal A$ are pairwise almost disjoint and $\mathcal A$ is not
a proper subset of a family $\mathcal A'$ with these properties.

Our main result is the following:

\textbf{Theorem 2: }Let $V=\underset{n<\omega}{\bigoplus} \mathbb{F}_2$
be a vector space of $\mathbb{F}_2$, then $V$ has no analytic mad
family of subspaces.

The rest of the paper will be devoted to the proof of Theorem 2. 

\textbf{Notation 2A: }a. $(x_n : n<\omega)$ will denote the basis
elements of $V$.

b. Given $u\subseteq \omega$, $\underset{n\in u}{\bigoplus} \mathbb{F}_2x_n$
will denote the subspace generated by $\{x_n : n\in u\}$.

c. For $v\in V$, the minimal $u\subseteq \omega$ such that $v\in \underset{n\in u}{\bigoplus} \mathbb{F}_2x_n$
will be denoted $supp(v)$.

d. For $u\subseteq \omega$, the subspace of $V$ generated by $\{ x_n : n\in u\}$
wil be denoted $span(u)$.

\textbf{Definition 3: }Given an ultrafilter $D$ on $V$, we define
the forcing $\mathbb{Q}=\mathbb{Q}_D$ as follows:

A. $p\in \mathbb Q$ iff $p=(u_p, \mathcal{A}_p)=(u, \mathcal A)$
where:

a. $u\subseteq V$ is finite and $0\in u$.

b. $\mathcal A \subseteq D$ is finite.

c. If $x\neq y \in u$ then the convex hulls of $supp(x)$ and $supp(y)$
are disjoint.

B. $(u_1,\mathcal{A}_1)\leq (u_2, \mathcal{A}_2)$ iff

a. $u_1 \subseteq u_2$ and for every $x\in u_1$ and $y\in u_2 \setminus u_1$,
$max(supp(x))<min(supp(y))$.

b. $\mathcal{A}_1 \subseteq \mathcal{A}_2$.

c. If $A \in \mathcal{A}_2 \setminus \mathcal{A}_1$ then 

$\alpha$. If $B\in \mathcal{A}_1$ then $A\subseteq B$.

$\beta$. If $B\in \mathcal{A}_1$, $x\in span(u_1)$ and $x+B \in D$
then $A \subseteq x+B$.

d. If $x\in span(u_1)$, $y\in span(u_2 \setminus u_1)$, $A\in \mathcal{A}_1$
and $x+A \in D$ then 

$\alpha$. $y\in (x+A) \cup \{0\}$.

$\beta$. $x+y+A \in D$.

\textbf{Definition 4: }For $D$ and $\mathbb Q=\mathbb{Q}_D$ as in
Definition 3, let $\underset{\sim}{A}=\cup \{ u_p : p\in \underset{\sim}{G}_{\mathbb Q} \}$.
If $|\underset{\sim}{A}|=\aleph_0$, we let $(\underset{\sim}{y_n} : n<\omega)$
be an enumeration of $\underset{\sim}{A}$ such that $max(supp(\underset{\sim}{y_n}))<min(supp(\underset{\sim}{}y_{n+1}))$.

\textbf{Observation 5: }$\mathbb Q$ is a partial order.

\textbf{Proof: }Suppose that $p_1 \leq p_2$ and $p_2 \leq p_3$,
we shall prove tht $p_1 \leq p_3$. Denote $p_l$ by $(u_l, \mathcal{A}_l)$
for $l=1,2,3$. Clauses (a) and (b) of Definition 3(B) are immediate.
Clause (c)$(\alpha)$ follows from (c)$(\beta)$, as $0\in span(u_1)$.
For (c)$(\beta)$, suppose that $A\in \mathcal{A}_3 \setminus \mathcal{A}_1$,
$B\in \mathcal{A}_1$, $x\in span(u_1)$ and $x+B \in D$. If $A \in \mathcal{A}_2 \setminus \mathcal{A}_1$,
then the desired conclusion follows from the fact that $p_1 \leq p_2$.
If $A \in \mathcal{A}_3 \setminus \mathcal{A}_2$, the result follows
similarly from the fact that $p_2 \leq p_3$. For clause (d), suppose
that $x\in span(u_1)$, $y\in span(u_3 \setminus u_1)$, $A\in \mathcal{A}_1$
and $x+A \in D$. As $y\in span((u_3 \setminus u_2) \cup (u_2 \setminus u_1))$,
there are $y_2 \in span(u_2 \setminus u_1)$ and $y_3 \in span(u_3 \setminus u_2)$
such that $y=y_2+y_3$. WLOG $y_3 \neq 0$, otherwise clause (d) follows
from the fact that $p_1 \leq p_2$. As $p_1 \leq p_2$, it follows
that $y_2 \in (x+A) \cup \{0\}$ and $x+y_2+A \in D$. As $p_2 \leq p_3$,
replacing $(x,y,A)$ in clause (d) by $(x+y_2,y_3,A)$ here, we get
that $y_3 \in x+y_2+A$ and $x+y_2+y_3+A \in D$. It follows that
$y=y_2+y_3 \in y_2+(x+y_2+A)=x+A$ (recall that $v+v=0$ for every
$v\in V$) and that $x+y+A \in D$. Therefore, $p_1 \leq p_3$. $\square$

\textbf{Observation 6: }If $p\in \mathbb Q$ and $B_1 \in D$, there
there is $q \in \mathbb Q$ and $B_2 \in \mathcal{A}_q$ such that
$p\leq q$ and $B_2 \subseteq B_1$.

\textbf{Proof: }Let $p\in \mathbb Q$ and $B_1 \in D$. Let $B_2:= \cap \{x+A : x\in span(u_p), A\in \mathcal{A}_p, x+A \in D \} \cap B_1$,
then $B_2 \in D$ and $B_2 \subseteq B_1$. Let $q=(u_p, \mathcal{A}_p \cup \{B_2\})$,
then $q\in \mathbb Q$ and it's easy to verify that $p\leq q$. $\square$

\textbf{Observation 7: }If $B\in D$ and every $p\in \mathbb Q$ can
be extended to $q \in \mathbb Q$ such that $|u_p|<|u_q|$ (so $\Vdash_{\mathbb Q} "|\underset{\sim}{A}|=\aleph_0"$)
then $\Vdash_{\mathbb{Q}} "$for some $k$ we have $span\{ \underset{\sim}{y_n} : k\leq n\} \subseteq B \cup \{0\}"$.

\textbf{Proof: }By the previous observation, there is a dense set
$I$ of conditions $p\in \mathbb Q$ with some $A\in \mathcal{A}_p$
such that $A\subseteq B$. Let $p\in I$, fix $A\in \mathcal{A}_p$
such that $A\subseteq B$ and let $k=|u_p|$, then $p$ forces values
$y_0,...,y_{k-1}$ to $\underset{\sim}{y_0},...,\underset{\sim}{y_{k-1}}$.
Suppose that $q$ is a condition above $p$ and let $y\in span(u_q \setminus u_p)$,
then $y\in A \cup \{0\} \subseteq B\cup \{0\}$ by Definition 3(B)(d)($\alpha$).
It follows that $p\Vdash_{\mathbb Q} "span\{ \underset{\sim}{y_n} : k\leq n\} \subseteq B \cup \{0\}"$.
As the last claim is true for any $p$ in the dense set $I$, this
completes the proof. $\square$

Towards the proof of Theorem 2, suppose that the theorem fails and
fix an analytic $V-$mad family $\mathcal A$ (i.e. $\mathcal A$
is a definition). We shall now derive a contradiction.

\textbf{Observation 8: }For $\mathbb Q$ as before,\textbf{ }$\Vdash_{\mathbb Q} "\mathcal A$
is $V-$mad family$"$.

\textbf{Proof: }Observe that as $V=\underset{n<\omega}{\bigoplus} \mathbb{F}_2$,
given two subspaces $S_1,S_2 \subseteq V$, $dim(S_1 \cap S_2)<\aleph_0$
iff $|S_1 \cap S_2|<\aleph_0$. As $\mathcal A$ is $\Sigma^1_1$,
the statement that $\mathcal A$ is maximal is $\Pi^1_2$ hence absolute.
Similarly, the almost disjointness of $\mathcal A$ is $\Pi^1_1$
and hence absolute. It follows that $\mathcal A$ is mad in $\bold{V}^{\mathbb Q}$.
$\square$

We shall now work with a forcing $\mathbb{Q}_D$ where $D$ is a certain
idempotent ultrafilter whose existence will be proved later.

\textbf{Fact 9: }There exists an ultrafilter $D$ on $V$ such that:

a. $D\cap \mathcal A=\emptyset$.

b. For every $A\in D$, for $D-$almost all $v$, $v+A\in D$.

c. $S\in D$ for every subspace $S$ of finite codimension.

Throughout the rest of the paper, $\mathbb Q=\mathbb{Q}_D$ where
$D$ is a fixed ultrafilter as in Fact 9 (which will be proved in
the end of the paper).

\textbf{Definition/Observation 10: }a. Let $\phi(A)=(\exists x)\psi(x, A)$
be the $\Sigma^1_1$ formula that defines $\mathcal A$. By the maximality
of $\mathcal A$ in $\bold{V}^{\mathbb Q}$, there are $\mathbb Q$-names
$\underset{\sim}{r}$ and $\underset{\sim}{B}$ such that $\Vdash_{\mathbb Q} " \psi(\underset{\sim}{r}, \underset{\sim}{B})$
and $|\underset{\sim}{B} \cap span\{ \underset{\sim}{y_n} : n<\omega\}|=\aleph_0"$
(by Observation 13 below, $\underset{\sim}{A}$ is infinite, hence
$\{\underset{\sim}{y_n} : n<\omega\}$ is an infinite well-defined
set).

b. $\Vdash_{\mathbb Q} "\underset{\sim}{B} \notin \bold{V}"$.

\textbf{Proof: }Let $A\in \mathcal{A}^{\bold{V}}$, then $A\notin D$
hence $V\setminus A \in D$. $\mathbb Q$ forces that $span\{ \underset{\sim}{y_n} : n<\omega\} \subseteq^* V\setminus A$.
As $\mathbb Q$ forces that $\underset{\sim}{B}$ contains infinitely
many elements from $span\{ \underset{\sim}{y_n} : n<\omega\}$, each
such element (modulo a finite number) is in $V\setminus A$ and it
follows that $\underset{\sim}{B} \neq A$. As $\underset{\sim}{B} \in \mathcal{A}^{\bold{V}^{\mathbb Q}}$,
it follows that $\underset{\sim}{B} \notin V$. $\square$

Let $\kappa=(2^{\aleph_0})^+$ and let $N$ be a countable elementary
submodel of $(H(\kappa),\in)$ such that $V,\phi, \underset{\sim}{r},\underset{\sim}{B} \in N$.
Let $(I_n : n<\omega)$ be an enumeration of the dense subsets of
$\mathbb Q$ that belong to $N$.

\textbf{Observation 11: }If $p\in \mathbb Q$ then $Z_p^+ \in D$
where $Z_p^+=\{v\in V : $ some $q\in \mathbb Q$ above $p$ forces
that $v\in \underset{\sim}{B} \cup \{0\}$, and moreover, $v\in span(u_q \setminus u_p)\}$.

\textbf{Proof: }Suppose towards contradiction that $Z_p^+ \notin D$,
then $V\setminus Z_p^+ \in D$. By the proof of Observation 6, there
is a condition $p_1$ above $p$ of the form $p_1=(u_p, \mathcal{A}_p \cup \{Z\})$
where $Z\subseteq V\setminus Z_p^+$. If $p_2$ is a condition above
$p_1$ such that $u_{p_2} \neq u_{p_1}$ (such a condition exists
by Observation 13), then $span(u_{p_2} \setminus u_p) \setminus \{0\} \subseteq Z\subseteq V\setminus Z_p^+$:
Let $y\in span(u_{p_2} \setminus u_p) \setminus \{0\}$, hence $y\in span(u_{p_2} \setminus u_{p_1})$.
As $Z\in \mathcal{A}_{p_1}$ and $p_1 \leq p_2$, it follows that
$y\in Z$, by the definition of the partial order. As $\Vdash_{\mathbb Q} "|\underset{\sim}{B} \cap span\{ \underset{\sim}{y_k} : k<\omega\}|=\aleph_0"$,
given $\{z_n : n<\omega\} \subseteq \underset{\sim}{B} \cap span\{ \underset{\sim}{y_k} : k<\omega\}$,
letting $z_n=z_n^1+z_n^2$ where $z_n^1 \in span\{ \underset{\sim}{y_k} : k\leq |u_{p_1}|\}$
and $z_n^2 \in span\{ \underset{\sim}{y_k} : k>|u_{p_1}| \}$, there
is an infinite set $\{n_k : k<\omega \}$ such that $z_{n_k}^1=z_*^1$
for all $k<\omega$. Therefore, for $k\neq k'$, as $\underset{\sim}{B}$
is a subspace, $z_{n_k}+z_{n_k'} \in \underset{\sim}{B}$ and moreover
$z_{n_k}+z_{n_k'} \in span \{ \underset{\sim}{y_k} : |u_{p_1}|<k\}$.
Therefore, there is $p_2$ above $p_1$ and $y$ such that $p_2 \Vdash_{\mathbb Q} "y\in span \{ \underset{\sim}{y_k} : k>|u_{p_1}|\} \cap \underset{\sim}{B} \setminus \{0\}"$.
By the definition of $\{ \underset{\sim}{y_k} : k<\omega\}$ and $\leq_{\mathbb Q}$,
it follows that $y\in span(u_{p_2} \setminus u_{p_1})=span(u_{p_2} \setminus u_p)$.
By the definition of $Z_p^+$, it follows that $y\in Z_p^+$. As $y\in span (u_{p_2} \setminus u_p) \setminus \{0\} \subseteq V\setminus Z_p^+$,
we obtain a contradiction. Therefore, $Z_p^+ \in D$. $\square$

\textbf{Observation 12: }If $p\in \mathbb Q$ then $Z_p^- \in D$
where $Z_p^-=\{ v\in V : $ some $q$ above $p$ forces $v\notin \underset{\sim}{B} \}$.

\textbf{Proof: }Suppose towards contradiction that $Z_p^- \notin D$,
then $\{v \in V : p \Vdash_{\mathbb Q} "v\in \underset{\sim}{B}"\}=V\setminus Z_p^- \in D$.
By the madness of $\mathcal{A}^{\bold{V}}$ in $\bold{V}$, there
is $B_1 \in \mathcal{A}^{\bold{V}}$ such that $B_2:=B_1 \cap (V\setminus Z_p^-)$
is infinite. Note that $p\Vdash "B_2 \subseteq V\setminus Z_p^- \subseteq \underset{\sim}{B}"$,
hence $p\Vdash_{\mathbb Q} "|\underset{\sim}{B} \cap B_1|=\aleph_0"$.
By absoluteness, $p\Vdash_{\mathbb Q} "B_1 \in \mathcal{A}"$, and
by the choice of $\underset{\sim}{B}$, $p \Vdash_{\mathbb Q} "\underset{\sim}{B} \in \mathcal A"$.
As $B_1 \in \bold{V}$ and $\Vdash_{\mathbb Q} "\underset{\sim}{B} \notin \bold{V}"$
it follows that $p\Vdash_{\mathbb Q} "\underset{\sim}{B} \neq B_1"$.
This contradicts the fact that $p\Vdash_{\mathbb Q} "\mathcal{A}$
is almost disjoint$"$. It follows that $Z_p^- \in D$, as required.
$\square$

\textbf{Observation 13: }$\Vdash_{\mathbb Q} "\underset{\sim}{A}$
is infinite$"$. Moreover, for every $p\in \mathbb Q$ there exists
$q\in \mathbb Q$ such that $p\leq q$ and $|u_p|<|u_q|$.

\textbf{Proof: }Let $p\in \mathbb Q$, we shall prove that there exists
$q\in \mathbb Q$ above $p$ such that $u_p \neq u_q$. Let $B_1=\cap \{ x+A : x\in u_p$,
$A \in \mathcal{A}_p$ and $x+A \in D\}$. As $D$ is a filter and
$u_p, \mathcal{A}_p$ are finite, $B_1 \in D$. By Fact 9, the set
$B_2=\{v \in V : v+B_1 \in D\}$ is in $D$. Let $n_*$ be large enough
such that $\underset{l<n_*}{\bigoplus} \mathbb{F}_2x_l$ includes
$u_p$, then by Fact 9, $B_3:=\underset{n_*<l}{\bigoplus} \mathbb{F}_2x_l \in D$.
Therefore, $B_1 \cap B_2 \cap B_3 \in D$ and hence is non-empty.
Let $y\in B_1 \cap B_2 \cap B_3$ and let $q=(u_p \cup \{ y\},\mathcal{A}_p)$.
Obviously, $q\in \mathbb Q$. It's easy to verify that $p\leq q$,
for example, we shall verify clause (B)(d)($\beta$) in Definition
3: Suppose that $x\in span(u_p)$, $A\in \mathcal{A}_p$ and $x+A \in D$.
$y\in B_2$, hence $y+B_1 \in D$. $B_1 \subseteq x+A$, hence $y+B_1 \subseteq y+x+A$.
As $y+B_1 \in D$, it follows that $y+x+A \in D$, as required. $\square$

\textbf{Finishing the proof of Theorem 2: }We shall now choose $(v_n, B_n^1, B_n^2, p_n^1, p_n^2)$
by induction on $n$ such that:

a. $v_n \in V$.

b. $B_n^l \subseteq V$ $(l=1,2)$.

c. $p_n^l \in N\cap \mathbb Q$ $(l=1,2)$.

d. For every $n<\omega$, $p_n^l \leq p_{n+1}^l$ $(l=1,2)$.

\textbf{Case I ($n=4i$): }We choose $p_{n+1}^l \in I_i$ above $p_n^l$
(recall that $I_i$ is dense).

\textbf{Case II ($n=4i+1$): }Suppose that $v_0,...,v_{n-1}$ and
$p_0^l,...,p_n^l$ have already been chosen. By Observation 11, $Z_{p_n^l}^+ \setminus \{v_0,...,v_{n-1}\} \in D$,
hence $(Z_{p_n^1}^+ \cap Z_{p_n^2}^+) \setminus \{v_0,...,v_{n-1}\} \in D$
and hence $(Z_{p_n^1}^+ \cap Z_{p_n^2}^+) \setminus \{v_0,...,v_{n-1}\} \neq \emptyset$.
Choose $v_n \in Z_{p_n^1}^+ \cap Z_{p_n^2}^+ \setminus \{ v_0,...,v_{n-1}\}$,
then there are conditions $q_{n+1}^1$ and $q_{n+1}^2$ above $p_n^1$
and $p_n^2$, respectively, such that $q_{n+1}^l \Vdash_{\mathbb Q} "v_n \in \underset{\sim}{B}"$
$(l=1,2)$. Let $p_{n+1}^l$ be an extension of $q_{n+1}^l$ that
decides $\underset{\sim}{B} \cap \underset{k<n}{\bigoplus}\mathbb{F}_2x_k$
and let $B_n^l$ be a subset of $V$ such that $p_{n+1}^l \Vdash_{\mathbb Q} "B_n^l=\underset{\sim}{B} \cap \underset{k<n}{\bigoplus}\mathbb{F}_2x_k"$.

\textbf{Case III ($n=4i+2$): }As in the previous case (using $Z_{p_n^1}^+$
and $Z_{p_n^2}^-$), we choose $v_n \notin \{v_0,...,v_{n-1}\}$ and
conditions $p_{n+1}^1, p_{n+1}^2$ such that $p_{n+1}^1 \Vdash_{\mathbb Q} "v_n \in \underset{\sim}{B}"$
and $p_{n+1}^2 \Vdash_{\mathbb Q} "v_n \notin \underset{\sim}{B}"$.

\textbf{Case IV ($n=4i+3$): }As in Case III (this time using $Z_{p_n^1}^-$
and $Z_{p_n^2}^+$), we choose $v_n \notin \{v_0,...,v_{n-1}\}$ and
conditions $p_{n+1}^1, p_{n+1}^2$ such that $p_{n+1}^1 \Vdash_{\mathbb Q} "v_n \notin \underset{\sim}{B}"$
and $p_{n+1}^2 \Vdash_{\mathbb Q} "v_n \in \underset{\sim}{B}"$.

Finally, having carried the induction, let $G_l=\{p\in N\cap \mathbb Q : p$
is below some $p_n^l\}$ $(l=1,2)$, then by Case I of the induction,
$G_l$ is $(N,\mathbb Q)$-generic. For $l=1,2$, let $S^l=\underset{n<\omega}{\cup}B_n^l$,
then by the genericity of $G_l$, the choice of the $B_n^l$s, and
$\mathcal A$ being analytic, it follows that $S_l \in \mathcal A$.
By Cases II-IV of the induction, $S_1 \neq S_2$ and $|S_1 \cap S_2|=\aleph_0$,
contradicting the almost disjointness of $\mathcal A$. This proves
Theorem 2 modulo Fact 9 that will be proved below. $\square$

\textbf{Proof of Fact 9: }For $S\in \mathcal{A}$ and $n<\omega$,
let $S[n]=\underset{l<n}{\bigoplus} \mathbb{F}_2 x_l +S$. Let $\mathcal D$
be the set of all nonprincipal ultrafilters on $V$ that contain all
subspaces of $V$ of finite codimension and all sets of the form $V\setminus S[n]$
for $S\in \mathcal A$ and $n<\omega$. Let $X=\{V\setminus S[n] : S\in\mathcal A, n<\omega\} \cup \{W : W\subseteq V$,
is a subspace of finite codimension$\}$, by the definition of the
topology on $\beta(V)$ (the space of ultrafilters on $V$), $\mathcal D$
is closed in $\beta(V)$. 

\textbf{Subclaim 1}: $\mathcal D \neq \emptyset$.

Proof: In order to show that $\mathcal D \neq \emptyset$, we shall
prove that $X$ has the FIP. Let $W\subseteq V$ be a subspace of
finite codimension, let $S_1,...,S_k \in \mathcal A$ and let $n_1,...,n_k<\omega$.
\textbf{Subclaim} \textbf{1(a)}: Given $S^1\neq S^2 \in \mathcal A$
and $n<\omega$, $S^1[n] \cap S^2[n]$ is finite. 

\textbf{Proof: }Let $k=|S^1 \cap S^2|$, we shall prove that $|S^1[n] \cap S^2[n]| \leq k2^{2n}$.
Suppose towards contradiction that $|S^1[n] \cap S^2[n]|>k2^{2n}$
and let $m:=k2^{2n}$. Let $\{r_j : j\leq m\}$ be pairwise distinct
elements of $S^1[n] \cap S^2[n]$. For each $j\leq m$ and $l\in \{1,2\}$,
there are $t_j^l$ and $a_{j,i}^l$ $(i<n)$ such that:

a. $t_j^l \in S^l$, $a_{j,i}^l \in \mathbb{F}_2$.

b. $r_j=\underset{i<n}{\Sigma}a_{j,i}^lx_i+t_j^l$.

Let $E$ be the equivalence relation on $\{j: j\leq m\}$ defined
by $j_1Ej_2$ iff $\underset{l=1,2}{\wedge} \underset{i<n}{\wedge}a_{j_1,i}^l=a_{j_2,i}^l$.
$E$ has $\leq 2^{2n}$ equivalence classes, hence there is $j_* \leq m$
such that $\frac{m+1}{2^{2n}} \leq |j_* /E|$, hence $k<|j_*/E|$.
By renaming, we may assume wlog that $\{0,1,....,k\} \subseteq j_*/E$.
For $l\in \{1,2\}$ and $j<k+1$, $r_j-r_0=t^l_j-t^l_0$, and as $t_j^l, t_0^l \in S^l$,
it follows that $r_j-r_0 \in S^l$. Therefore, $r_j-r_0 \in S^1 \cap S^2$
for every $j<k+1$. As $\{r_j : j<k+1\}$ is without repetition, so
is $\{r_j-r_0 : j<k+1\}$, contradicting the fact that $|S^1 \cap S^2|=k$.
This proves Subclaim 1(a).

Back to the proof of Subclaim 1, choosing $S' \in \mathcal A \setminus \{S_1,...,S_k\}$,
$S' \cap S_l[n_l] \subseteq S'[n_l] \cap S_l[n_l]$ is finite. WLOG
suppose that $W=\underset{m\leq n}{\bigoplus} \mathbb{F}_2x_n$ and
let $\{z_n : n<\omega\}$ be an infinite subset of $S'$ such that
$max(supp(z_n))<max(supp(z_{n+1}))$, wlog $\{z_n : n<\omega\}$ is
disjoint to $S_l[n_l]$ for every $l\leq k$. By the same argument
as in the proof of Observation 11, we may assume wlog that for each
$n$, $z_n=z \oplus z_n'$ for a fixed $z\in \underset{l<m}{\bigoplus} \mathbb{F}_2x_l$.
Now consider the set $\{z_0+z_n : n<\omega\} \subseteq S'$. As $max(supp(z_n))<max(supp(z_{n+1}))$,
this set is infinite, and therefore contains an element $z_0+z_n$
such that $z_0+z_n \notin \underset{l \leq k}{\cup}S_l[n_l]$. Obviously,
$z_0+z_n \in W$. Therefore, $W\cap (\underset{l\leq k}{\cap}(V\setminus S_l[n_l])) \neq \emptyset$
and it follows that $X$ has the FIP. This completes the proof of
Subclaim 1.

\textbf{Subclaim 2:} If $D_1,D_2 \in \mathcal D$ then $D_1 \oplus D_2 \in \mathcal D$
where $D_1 \oplus D_2=\{A : (\forall^{D_1}s_1)(\forall^{D_2}s_2)(s_1+s_2 \in A)\}$.

\textbf{Proof: }Obviously, $D_1 \oplus D_2 \in \beta(V)$. For every
space of the form $W=\underset{m\leq n}{\bigoplus}\mathbb{F}_2x_n$
we have $W\in D_1$, $W\in D_2$ and $s_1+s_2 \in W$ for every $s_1,s_2 \in W$.
Therefore, $W\in D_1 \oplus D_2$ and the ultrafilter contains all
subspaces of finite codimension. Now let $S\in \mathcal A$ and $n<\omega$,
we shall prove that $S[n] \notin D_1 \oplus D_2$. Fix $s_1 \in V$,
letting $m>n$ such that $s_1 \in \underset{i<m}{\bigoplus} \mathbb{F}_2x_i$,
we have $\{ s_2 : s_1+s_2 \in S[n]\} \subseteq S[m]$ (if $s_1+s_2 \in S[n]$,
then $s_1+s_2=x+y$ for some $x\in \underset{i<n}{\bigoplus}\mathbb{F}_2x_i$
and $y\in S$. Hence, $s_2=s_1+s_1+s_2=(s_1+x)+y \in \underset{i<m}{\bigoplus} \mathbb{F}_2x_i + S=S[m]$).
As $S[m] \notin D_2$, it follows that $S[n] \notin D_1 \oplus D_2$.
This completes the proof of Subclaim 2.

By Glazer's argument in the proof of Hindman's theorem (see Lemma
10.1, page 449 in {[}Co{]} or Lemma 2.7 in {[}RoSh:957{]}), the fact
that $\mathcal D$ is a nonempty closed subset of $\beta(V)$ that
is closed under the $\oplus$-operation implies that there exists
$D\in \mathcal D$ such that $D\oplus D=D$. We shall now check that
$D$ is as required in Fact 9. Clauses (a) and (c) are immediate from
the fact that $D\in \mathcal D$, so it remains to show that $D$
satisfies clause (b). Let $A\in D$, we shall prove that $v+A \in D$
for $D$-almost all $v$. Given $v\in V$, let $v\oplus A=\{z \in V : v+z \in A\}$.
As $D\oplus D=D$, $(\forall^D s_1)(\forall^D s_2)(s_1+s_2 \in A)$.
Therefore, $(\forall^D s_1)(s_1 \oplus A \in D)$. Note that $s_1 \oplus A=\{s_2 : s_2 \in s_1+A\}=s_1+A$,
therefore, for $D$-almost all $s_1$ we have $s_1+A\in D$, as required.
This completes the proof of Fact 9 $\square$

The following question remains open:

\textbf{Question: }Let $F\neq \mathbb{F}_2$ be a field of cardinality
$\leq \aleph_0$, is there an analytic mad families of subspaces of
$\underset{n<\omega}{\bigoplus} F$?

It is conceivable that a method similar to that of {[}HwSh:1089{]}
might allow us to construct a Borel mad family for fields other than
$\mathbb{F}_2$.

Finally, we observe that combining the proof of {[}HwSh:1090{]} with
the results from this paper we obtain the following:

\textbf{Theorem: }Let $V=\underset{n<\omega}{\bigoplus} \mathbb{F}_2$,
then $ZF+DC+"$there are no mad families of subspaces of $V"$ is
equiconsistent with $ZFC$.

The proof is almost identical to {[}HwSh:1090{]}, where instead of
using Mathias forcing, we now use the forcing $\mathbb{Q}_D$ from
this paper where $D$ is as in Fact 9. We shall elaborate on the proof
in a subsequent paper.

\textbf{\large References}{\large \par}

{[}Co{]} W. Wistar Comfort, Ultrafilters: some old and some new results.
Bulletin of the American Mathematical Society, 83:417-455, 1977

{[}HwSh:1089{]} Haim Horowitz and Saharon Shelah, A Borel maximal
eventually different family, arXiv:1605.07123

{[}HwSh:1090{]} Haim Horowitz and Saharon Shelah, Can you take Toernquist's
inaccessible away? arXiv:1605.02419

{[}HwSh:1095{]} Haim Horowitz and Saharon Shelah, A Borel maximal
cofinitary group, arXiv:1610.01344

{[}Ma{]} A. R. D. Mathias, Happy families, Ann. Math. Logic \textbf{12
}(1977), no. 1, 59-111

{[}RoSh:957{]} Andrzej Roslanowski and Saharon Shelah, Partition theorems
from creatures and idempotent ultrafilters, Annals Combinatorics 17
(2013) 353-378

{[}Sm{]} Iian Smythe, Madness in vector spaces, arXiv:1712.00057

$\\$

(Haim Horowitz) Department of Mathematics

University of Toronto

Bahen Centre, 40 St. George St., Room 6290

Toronto, Ontario, Canada M5S 2E4

E-mail address: haim@math.toronto.edu

$\\$

(Saharon Shelah) Einstein Institute of Mathematics

Edmond J. Safra Campus,

The Hebrew University of Jerusalem.

Givat Ram, Jerusalem 91904, Israel.

Department of Mathematics

Hill Center - Busch Campus,

Rutgers, The State University of New Jersey.

110 Frelinghuysen Road, Piscataway, NJ 08854-8019 USA

E-mail address: shelah@math.huji.ac.il
\end{document}